\documentclass[12pt,a4paper]{article}
\usepackage[intlimits]{amsmath}
\usepackage{amsfonts,amssymb,amscd,amsthm}

\theoremstyle{plain}
\newtheorem{theorem}{Theorem}

\newtheorem{lemma}[theorem]{Lemma}

\theoremstyle{definition}

\theoremstyle{remark}
\newtheorem{remark}[theorem]{Remark}

\def\RR{\mathbb{R}}

\def\bcdot{\,\boldsymbol\cdot\,}

\def\ga{\gamma}
\def\la{\lambda}
\def\si{\sigma}

\def\Ga{\Gamma}
\def\La{\Lambda}

\DeclareMathOperator\dist{dist}

\DeclareMathOperator\tr{tr} \DeclareMathOperator\ord{ord}
\DeclareMathOperator\ind{ind}

\newcommand{\op}[1]{\operatorname{#1}}

\def\pa{\partial}

\let\rom\textup

\newcommand{\BL}{\biggl}
\newcommand{\BR}{\biggr}

\newcommand{\RRe}{{\rm Re}}
\newcommand{\abs}[1]{\lvert#1\rvert}
\newcommand{\norm}[1]{\left\| #1 \right\|}

\def\Td{\operatorname{Td}}
\def\Pf{\operatorname{Pf}}
\def\ch{\operatorname{ch}}

\title{On Elliptic Differential Operators\\[2pt] with Shifts}
\author{V.~E.~Nazaikinskii, A.~Yu.~Savin, and B.~Yu.~Sternin}
\date{}

\begin{document}
\maketitle

\section{Statement of the Problem \\ and Main Theorem}

The motivating point of our research was differential
operators on the non\-com\-mu\-ta\-ti\-ve
torus\footnote{Which have important applications to
specific physical problems such as quantum Hall effect.}
studied by Connes~\cite{MR572645,Con1}, who in particular
obtained an index formula for such operators. These
operators include shifts (more precisely, in this case,
irrational rotations); hence our interest in general
differential equations with shifts naturally arose.

Let $M$ be a smooth closed manifold. We consider operators
on $M$ of the form
\begin{equation}\label{eq1}
    Du =\sum_j(g_j^*)^{-1}(D_ju),
\end{equation}
where the sum is finite, $D_j$ are differential operators, and
$g_j^*$ are shift operators,
\begin{equation}
    g_j^*u(x)=u(g_j(x)).
\end{equation}
Here $g_j\colon M\to M$ are some diffeomorphisms. Elliptic
theory for operators of the form~\eqref{eq1} is well known
(e.g., see \cite{AnLe3}). In particular, under certain
assumptions about $g_j$, the principal symbol of the
operator~\eqref{eq1} is defined and its invertibility
implies that the operator~\eqref{eq1} is Fredholm in
appropriate Sobolev spaces on~$M$. However, apart from
Connes' example, the index formula for such operators has
so far been obtained only for the case in which the group
$\Ga$ generated by the elements $g_j$ is \textit{finite},
and neither the formula itself nor the proof method make
sense for $\Ga$ infinite. However, it is the case of an
\textit{infinite} group that is of main interest in
applications, as, for example, in Connes' above-mentioned
work. In the present paper, we prove an index formula for
the case of an \textit{infinite} group $\Ga$ satisfying
certain mild assumptions.

\paragraph{Acknowledgements.}
The research was supported in part by RFBR grants
nos.~05-01-00982 and~06-01-00098 and DFG grant 436 RUS
113/849/0-1\circledR ``$K$-theory and noncommutative
geometry of stratified manifolds."

The authors thank Professor Schrohe and Leibniz
Universit\"at Hannover for kind hospitality.

\subsection{Definition of the operator algebra}

\subsubsection{The group $\Ga$}
We assume that $M$ is a compact oriented Riemannian
manifold and $\Ga$ is a finitely generated dense subgroup
of a compact Lie group of orientation-preserving isometries
of $M$ satisfying the following two conditions:
\begin{enumerate}
    \item (The \textit{polynomial growth} condition.) For some
    finite system of generators of $\Ga$, the number of distinct
    elements of $\Ga$ representable by words of length $\le k$ in
    these generators grows not faster than some power of $k$.

    In what follows, we fix some system of generators and denote by
    $\abs{g}$ the minimum length of words representing $g\in\Ga$.
    \item (The ``\textit{Diophantine}'' condition.) There exists a
    finite $N$ and a constant $C>0$ such that the estimate
\begin{equation*}
    \dist(g(x),x)\ge C\abs{g}^{-N}\dist(x,\op{fix}(g))
\end{equation*}
holds for all $x\in M$ and $g\in\Ga$. Here
$\dist(\bcdot,\bcdot)$ stands for the Riemannian distance
between points of the manifold, and $\op{fix}(g)$ is the
set of fixed points of $g\in\Ga$. If $\op{fix}(g)$ is
empty, then by convention we set $\dist(x,\op{fix}(g))=1$.
\end{enumerate}

Each element $g\in\Ga$ naturally acts on the cotangent space $T^*M$
(it takes each fiber $T_x^*M$ to the fiber $T_{g(x)}^*$ by the
linear transformation $((dg)^*)^{-1}$, where $dg$ is the
differential of $g$); this action will be denoted by $\pa g$. Since
$g$ is an isometry, it follows that $\pa g$ is also an isometry and
in particular has a well-defined restriction to $S^*M$.

\subsubsection{Pseudodifferential operators}

Let $D_g$, $g\in\Ga$, be a family of pseudodifferential
operator of order $\le m$ on $M$ rapidly (i.e., faster than
any power of $\abs{g}^{-1}$) decaying as $\abs{g}\to\infty$
in the Fr\'echet topology on the set of pseudodifferential
operators (e.g., see \cite{EgSc1}). In particular, it
follows that
\begin{equation*}
    \norm{D_g\colon H^s(M)\to H^{s-m}(M)}\le C_{sN}(1+\abs{g})^{-N}
\end{equation*}
for any $s$ and $N$. Then the series
\begin{equation}\label{eq2}
    D=\sum_{g\in\Ga}(g^*)^{-1}\circ D_g\colon H^s(M)\to H^{s-m}(M)
\end{equation}
converges absolutely in operator norm for every $s\in\RR$. The set
of such operators will be denoted by $\Psi_\Ga^m$, and we set
\begin{equation*}
 \Psi_\Ga^\infty=\bigcup_m\Psi_\Ga^m .
\end{equation*}
One can readily prove that under our assumptions
$\Psi_\Ga^\infty$ is an algebra. It is called the
\textit{algebra of $\Psi$DO with shifts} on $M$.

\subsubsection{Symbols and the Fredholm property}

For the operator~\eqref{eq2}, we define its (\textit{principal})
\textit{symbol} by the formula
\begin{equation}\label{sima1}
\sigma(D):=\sum_{g\in \Gamma}(\pa g^*)^{-1}\circ
\sigma(D_g) \colon L^2(S^*M)\longrightarrow L^2(S^*M),
\end{equation}
where $\si(D_g)$ is the conventional principal symbol of the
$m$th-order pseudodifferential operator $D_g$, acting as a
multiplication operator.

\begin{theorem}[\cite{AnLe3}]\label{th1}
The principal symbol is well defined, and the operator
$$
D:H^s(M)\longrightarrow H^{s-m}(M),
$$
where $D\in\Psi^m_\Gamma$, is Fredholm if and only if its symbol
$\sigma(D)$ is invertible. Its index is independent of $s$.
\end{theorem}

The expression~\eqref{sima1} belongs to the algebra
$C^\infty(S^*M)_\Gamma$ of $C^\infty(S^*M)$-valued
functions on $\Ga$ rapidly decaying in the Fr\'echet
topology of $C^\infty(S^*M)$ as $\abs{g}\to\infty$.
By~\cite{Schwe1,Schwe2}, this is a dense local subalgebra
in $C(S^*M)_\Gamma\equiv C(S^*M)\rtimes\Gamma$.
Consequently, the inverse of an elliptic symbol $\si\in
C^\infty(S^*M)_\Gamma$ also belongs to
$C^\infty(S^*M)_\Gamma$.

\subsection{Main Theorem}

\subsubsection{Notation}

To state the index theorem, we need some notation.

If $E$ is a vector bundle over a manifold, then by
$\lambda_{-1}(E)$ we denote the (virtual) vector bundle
\begin{equation*}
 \lambda_{-1}(E):=\sum_{j\ge 0}(-1)^j\Lambda^j(E)
 =\La^{even}(E)-\La^{odd}(E),
\end{equation*}
composed of the exterior powers of $E$ (the sum is actually
finite). By
\begin{equation*}
 \ch E (g):=\tr\BL( g^* \exp\BL(-\frac 1{2\pi
i}\Omega\BR)\BR)
\end{equation*}
we denote the localized Chern character of a bundle $E$
with some fixed curvature form $\Omega$ at an element
$g\in\Gamma$. Here $\tr$ stands for the trace in the fibers
of the vector bundle.

Next, let us define the Chern--Simons character of an
elliptic element $a\in C^\infty(S^*M)_\Gamma$. Consider the
algebra $\Lambda(S^*M)_\Gamma$ of $\Lambda(S^*M)$-valued
functions on $\Ga$ rapidly decaying in the Fr\'echet
topology of $\Lambda(S^*M)$ (with the obvious
multiplication). This algebra is graded (by form degree)
and becomes a differential graded algebra if we equip it
with the differential
\begin{equation}\label{difa1}
    d\BL(\sum_g(\pa g^*)^{-1}\circ\omega_g \BR)
    := \sum_g(\pa g^*)^{-1}\circ (d\omega_g).
\end{equation}
Now we can define the Chern--Simons character by setting
\begin{equation}\label{cs1}
\begin{aligned}
    \ch a&:=\ch_1 a+\ch_3 a+\ldots,\quad\text{where}\\
    \ch_{2k+1}a&:= \BL(\frac 1{2\pi i}\BL)^{k+1}
     \frac{(k)!}{(2k+1)!} \tr [a^{-1}da]^{2k+1}.
\end{aligned}
\end{equation}
By $\ch a(g)$ we denote the differential form that is the
coefficient of $g^*$ in the expansion
\begin{equation*}
    \ch a=\sum_{g\in\Ga}(\pa g^*)^{-1}\circ\ch a(g).
\end{equation*}

Next, recall that if $g$ is an isometry of $M$, then the
set $\op{fix}g$ of fixed points of $g$ is a disjoint union
of smooth submanifolds of $M$ \cite{AtSe2}, each of which
is locally (i.e., in a neighborhood of its every point) the
image, under the geodesic exponential mapping, of the
eigenspace of $dg$ corresponding to the eigenvalue~$1$.
These submanifolds will be denoted by $M_g$. (We shall not
use double subscripts to avoid clumsiness.) Our index
formula involves integration over the cosphere bundles
$S^*M_g$; we equip these with the standard orientation
coming from the almost complex structure on $T^*M_g$. By
$NM_g$ we denote the normal bundle of $M_g$.

\subsubsection{The index formula}

\begin{theorem}
Let $D\in \Psi^\infty(M,Mat(m,\mathbb{C}))_\Gamma$ be a matrix
elliptic operator on $M$. Then the following index formula
holds\rom:
\begin{equation}\label{fixedp}
\ind D=\sum_{g\in \Gamma}\BL[\sum_{M_g}\int_{S^*M_g}
\frac{\Td(T^*M_g\otimes
\mathbb{C})}{\ch\lambda_{-1}(NM_g\otimes
\mathbb{C})(g)}{\ch\sigma(D)(g)}\BR].
\end{equation}
Here the denominators do not vanish, and the double series
converges absolutely.
\end{theorem}

\begin{remark}
The nonvanishing of the denominators was proved in \cite{AtSe2}.
However, the convergence of the series is still to be proved.
\end{remark}

\section{Proof of the Index Formula}

Let $D\in \Psi^\infty(M,Mat(m,\mathbb{C}))_\Gamma$ be a matrix
elliptic operator on $M$. For simplicity, we assume that $D$ is a
first-order differential operator. The general case can be treated
by order reduction and by using the technique of operators with
continuous symbols, just as in the classical paper \cite{AtSi1}.

Our proof is based on the reduction of the original operator on $M$
to a Dirac type operator on $S^*M$ with the use of a trick that is
apparently due to Kasparov. (However, we have not been able to find
anything of this sort in the literature and so cannot give a
precise reference.)

\subsection{Reduction to an operator on the cosphere bundle}

We first perform reduction to the cotangent bundle and only then
use some $K$-theory to go down to the cosphere bundle.

\subsubsection{Elliptic theory on the cotangent bundle}

We consider symbols $\sigma(x,\xi,p,\eta)$, where $(x,\xi)$ are
standard coordinates on $T^*M$ and $(p,\eta)$ are the dual momenta,
satisfying the estimates
$$
|\partial^\alpha_x\partial^\beta_z \sigma|\le
C_{\alpha\beta}(1+|z|)^{m-|\beta|},
$$
where $z=(\xi,p,\eta)$. We also introduce associated
principal symbols, homogeneous in $z$. The corresponding
operators are considered in spaces with the norms
$$
   \|u\|^2_s:
   =\int_{T^*M}|\left(\Delta_x+\Delta_\xi+\xi^2\right)^{s/2}u|^2
        dx\,d\xi,
$$
where $\Delta_x$ and $\Delta_\xi$ are positive Laplace
operators in the corresponding variables. The ellipticity
condition is that the principal symbol should be invertible
for $z\ne 0$.

The group $\Gamma$ acts on $T^*M$ by the isometries $\pa g$. Hence
on can define, in a manner similar to the above, an algebra of
$\Psi$DO with shifts on $T^*M$, which will be denoted by
$\Psi(T^*M)_\Gamma$.

\subsubsection{Reduction to the cotangent bundle}

In the space $\mathbb{R}^n$, $n=\dim M$, consider the elliptic
operator
\begin{multline*}
E=\left[(x+\partial/\partial x)dx\right]
\wedge+\left[(x-\partial/\partial x)dx \right]\lrcorner\colon
C^\infty(\mathbb{R}^n,\Lambda^{ev}(\mathbb{C}^n))
\\\longrightarrow
C^\infty(\mathbb{R}^n,\Lambda^{odd}(\mathbb{C}^n)),
\end{multline*}
where
$$
xdx:=\sum_j x_j dx_j;\qquad \frac\partial{\partial x}dx:= \sum_j
\frac\partial{\partial x_j} dx_j.
$$
This operator is elliptic of index one in appropriate Sobolev
spaces; its kernel is one-dimensional and is spanned by the
function $\exp(-x^2/2)$. Moreover, the operator is
$O(n)$-invariant. Hence we can consider the family
$\mathcal{E}=\{\mathcal{E}_x\}_{x\in M}$ of such operators acting
in the fibers of the cotangent bundle,
$$
\mathcal{E}:C^\infty(T^*M,\pi^*\Lambda^{ev}(M))\longrightarrow
C^\infty(T^*M,\pi^*\Lambda^{odd}(M)),\quad \pi:T^*M\longrightarrow
M.
$$
In turn, the operator $D$ can be lifted to the operator
$$
D\otimes
1_\Lambda:C^\infty(T^*M,\mathbb{C}^m\otimes\pi^*\Lambda(M))\longrightarrow
C^\infty(T^*M,\mathbb{C}^m\otimes\pi^*\Lambda(M)).
$$

\begin{lemma}
The crossed product
\begin{multline*}
D\# \mathcal{E}:=\left(%
\begin{array}{cc}
  D\otimes 1_{\Lambda^{ev}} & -1_{\mathbb{C}^m}\otimes\mathcal{E}^* \\
  1_{\mathbb{C}^m}\otimes\mathcal{E} & D^*\otimes 1_{\Lambda^{odd}} \\
\end{array}%
\right): C^\infty(T^*M,\mathbb{C}^m\otimes
\pi^*\Lambda(M))\\
\longrightarrow C^\infty(T^*M,\mathbb{C}^m\otimes
\pi^*\Lambda(M))
\end{multline*}
is an elliptic operator in $\Psi(T^*M)_\Gamma$, and one has
$$
\ind D=\ind (D\# \mathcal{E}).
$$
\end{lemma}
\begin{proof}
The proof is by analogy with that of a similar assertion in
\cite{AtSi1}.
\end{proof}

\subsubsection{Reduction to the Todd operator}

We continue the symbol of $D$ as a first-order homogeneous
function to the entire space $T^*M$.

\begin{lemma}\label{ktheorem}
One has
\begin{equation}\label{dir1}
\ind (D\#\mathcal{E})=\ind (\sigma(D)\# \mathcal{D}),
\end{equation}
where $\mathcal{D}$ is the Todd operator on the space $T^*M$,
viewed as an almost complex manifold with complex coordinates
$$
z_1=\xi_1+ix_1,\;...,\; z_n=\xi_n+ix_n.
$$
\end{lemma}
\begin{proof}
We denote the variables on $T^*M$ by $(x,\xi)$ and the dual momenta
by $(p,\eta)$.

The symbols of the operators in~\eqref{dir1} belong to the algebra
$$
C^\infty(S(T^*M\oplus T^*M\oplus T^*M))_\Gamma,
$$
where $\Gamma$ acts on the cosphere bundle of $T^*M$ by the second
differential
$$
\pa^2g:=\pa(\pa g):T^*(T^*M)\longrightarrow T^*(T^*M).
$$
Let $D=\sum D_g g^*$. Then the symbols of the operators on the
left- and right-hand sides in \eqref{dir1} are
\begin{gather*}
\BL(\sum_g(\pa^2g^*)^{-1}\circ\sigma(D_g)(x,p)\BR)\#
\sigma(\mathcal{E})(\xi,\eta),\\
\BL(\sum_g(\pa^2g^*)^{-1}\circ\sigma(D_g)(x,\xi)\BR)\#
\sigma(\mathcal{E})(-p,\eta),
\end{gather*}
respectively. These symbols can be joined by the homotopy
$$
\BL(\sum_g(\pa^2g^*)^{-1}\circ\sigma(D_g)(x,\xi\sin\varphi
+p\cos\varphi )\BR)\# \sigma(\mathcal{E})(\xi \cos\varphi
 -p\sin\varphi,\eta),
$$
$\varphi\in[0,\pi/2]$, and one can readily verify that this
homotopy preserves ellipticity. The proof of the lemma is complete.
\end{proof}

\subsubsection{Passage to the cosphere bundle}

The Todd operator on $T^*M$ induces the so-called
\emph{tangential Todd operator}
$$
\mathcal{D}_{S^*M}: C^\infty(S^*M,\Lambda^{ev}(M))\longrightarrow
C^\infty(S^*M,\Lambda^{ev}(M))
$$
on the cosphere bundle $S^*M\subset T^*M$. This is an
elliptic symmetric operator. To describe its principal
symbol, we identify the tangent and cotangent bundle using
the metric and take local coordinates $(x,\xi,p,\eta)$ on
$T^*(T^*M)$. Let $c(p,\eta):=\sigma
(\mathcal{D})(x,\xi,p,\eta),$ where $\mathcal{D}$ is the
Todd operator; then the symbol of the tangential Todd
operator is given by
$$
\sigma(\mathcal{D}_{S^*M})(x,\xi,p,\eta):=ic(0,\xi)c(p,\eta),
$$
where $\xi$, $|\xi|=1$, is a unit covector and $(p,\eta)\in
T_{(x,\xi)}^*S^*M,$ $\eta\perp \xi$, is a covector tangent
to $S^*M$.

A standard $K$-theoretic argument shows that the following lemma is
true.
\begin{lemma}
One has
$$
\ind \left[\sigma(D)\# \mathcal{D}\right]=\ind  \left[P
\sigma(D)|_{S^*M}P\right],
$$
where
\begin{equation}\label{topa1}
P \sigma(D)|_{S^*M}P: PL^2(S^*M,\Lambda^{ev}(M)\otimes
\mathbb{C}^m) \longrightarrow PL^2(S^*M,\Lambda^{ev}(M)\otimes
\mathbb{C}^m)
\end{equation}
is a Toeplitz operator on the subspace determined by the positive
spectral projection $P$ of the tangential Todd operator on the
cosphere bundle.
\end{lemma}

Thus, to compute the index of the original operator $D$, it remains
to compute the index of the Toeplitz operator \eqref{topa1} on
$S^*M$.

\subsection{Computation of the index of a Dirac type operator}

The Todd operator that arose in the preceding subsection is
a Dirac type operator with noncommuting coefficients. In
this section, we show how to compute the index of such
operators. First, we deal with the case of even-dimensional
manifolds, and then use the standard technique of
multiplication by $\mathbb{S}^1$ to cover the case of
$S^*M$, which is actually of interest to us.

\subsubsection{The even case}

Let $X$ be a smooth closed manifold on which the group
$\Gamma$ acts isometrically. We assume that the power
growth and Diophantine conditions are satisfied and that
the manifold is even-dimensional and is equipped with a
spin structure, which is moreover preserved by the action
of $\Ga$. By
$$
D_+:C^\infty(X,S_+)\longrightarrow C^\infty(X,S_-)
$$
we denote the corresponding $\Gamma$-invariant Dirac operator on
$X$. This operator is local with respect to the action of the
algebra $C^\infty(X)_\Gamma$, and hence we obtain an elliptic
operator if we twist $D_+$ by some projection
$$
p\in Mat(N, C^\infty(X)_\Gamma)
$$
over the algebra $C^\infty(X)_\Gamma$ of rapidly decaying
$C^\infty(X)$-valued functions on the group.

We denote the twisted operator by
\begin{equation}
\label{twist1} pD_+p: pC^\infty(X,S_+\otimes
\mathbb{C}^N)\longrightarrow pC^\infty(X,S_-\otimes \mathbb{C}^N).
\end{equation}

To give the index formula, we introduce the noncommutative
Chern character of the projection $p$. Consider the graded
differential algebra $\Lambda(X)_\Gamma$. The
noncommutative Chern character is given by
$$
\ch p:=\tr p\BL[\exp\BL(-\frac 1{2\pi i}pdpdp\BR)\BR] \in
\Lambda(X)_\Gamma,
$$
where, as above, $\tr$ is the usual matrix trace.

\begin{theorem}
One has
\begin{equation}\label{dirac1}
 \ind pD_+p=\sum_{g\in \Gamma} \BL[\sum_{X_g\subset X}
\int_{X_g}
{A(X_{g})}\Pf\left\{2i\sin\left(\frac\Omega{4\pi}+\frac{i\Theta}2\right)(NX_{g})\right\}^{-1}
\cdot \ch p(g)\BR],
\end{equation}
where the double sum is absolutely convergent. Here
\begin{itemize}
    \item $X_g$ are connected components of the set of fixed points
          of $g$, and each $X_g$ is a smooth even-dimensional
          manifold\rom;
    \item $\Omega(NX_g)$ is the curvature form of the normal bundle
          $NX_g$\rom;
    \item $\Theta(NX_g)$ is the logarithm of Jacobi's matrix of
          the mapping $g^*:NX_g\to NX_g$\rom;
    \item $\Pf$ is the Pfaffian of a skew-symmetric matrix\rom;
    \item $\ch p(g)$ is the coefficient of $(g^*)^{-1}$ in the
          decomposition of the Chern character
          as an element of $\Lambda(X)_\Gamma$\rom;
    \item $A(X_g)$ is a differential form representing
          the $A$-class of the manifold $X_g$.
\end{itemize}
\end{theorem}

\begin{proof}
For simplicity, we assume that the Dirac operator $D_+$ is
invertible. The proof is based on the application of the
Connes--Moscovici formula for the computation of the
Chern--Connes character \cite{CoMo1} of the operator $D_+$
with the subsequent computation of the terms in this
formula with the use of the explicit formulas given in
\cite[Theorem 5, p.~471]{CheHu1}. These computations are
standard but clumsy, and we omit them. Instead, we show
that the Connes--Moscovici formula applies to our case.

Let $D:C^\infty(X,S_+\oplus S_-)\longrightarrow
C^\infty(X,S_+\oplus S_-)$ be the full Dirac operator
$$
D=\left(%
\begin{array}{cc}
  0 & (D_+)^* \\
  D_+ & 0 \\
\end{array}%
\right).
$$
By the results given in \cite{Hig6}, the Connes--Moscovici formula
applies if the following two conditions are satisfied for the
spectral triple $(C^\infty(X)_\Gamma,L^2(X,S),D)$ determined by the
Dirac operator:

1) It has \emph{finite analytic dimension}; i.e., for some
$d,q\in \mathbb{R}$ and all differential operators $P\in
\Psi_\Gamma^\infty$ the function $\tr P\Delta^{-z}$ is
holomorphic in the half-plane $Re z>(\ord P-q)/d$ (here
$\Delta:=D^*D$).

2) It has \emph{the analytic continuation property}; i.e.,
for each differential operator $P\in {\rm
Diff}(X)_\Gamma^\infty$ the function $\tr P\Delta^{-z}$
admits a meromorphic continuation into the entire complex
plane.

The first assumption obviously holds. (It is valid for
operators without shifts, and for operators with shifts the
proof is the same.)

Let us verify the second assumption. Let
$P=\sum_g(g^*)^{-1}\circ P_g$, where the coefficients $P_g$
rapidly decay. We should prove that

1. Each function $\tr (g^*)^{-1}\circ P_g \Delta^{-z}$
admits a meromorphic continuation into the entire complex
plane with the same discrete set of possible poles.

2. The series formed by these functions locally uniformly
converges outside this set, thus giving a meromorphic
function as the sum.

To this end, one uses the representation
\begin{equation*}
    \Delta^{-z}=\frac1{2\pi i}\int_\gamma
    \la^{-z}(\la-\Delta)^{-1}dz
    =\frac{l!}{(1-z)\dotsm(l-z)}\int_\gamma
    \la^{l-z}(\la-\Delta)^{-1-l}dz,
\end{equation*}
where $\ga$ is a contour going along the imaginary axis and
bypassing the origin on the right. For sufficiently large
$l$, the operator $(g^*)^{-1}\circ P_g (\la-\Delta)^{-1-l}$
is trace class, and the integral converges and defines an
analytic function for $\RRe z\gg0$. To get the desired
analytic continuation, we note that $(\la-\Delta)^{-1-l}$
is a pseudodifferential operator with parameter $\la$ in
the sense of \cite{Shu1} and hence has the complete symbol
possessing an asymptotic expansion in homogeneous functions
of integer orders of homogeneity in $(\xi,\la)$, where
$\xi$ are the covariables. In the trace integral, only an
arbitrarily small neighborhood of $\op{fix}g$ can produce
terms that are not extendable as entire functions; in these
neighborhoods, the stationary phase method reduces the
integrals involving the homogeneous components of the
symbol over the phase space to asymptotic expansions in
(half)-integer powers of $\la$ whose coefficients are
integrals, independent of $\la$, over the cotangent bundles
of $M_g$. The subsequent integration over $\gamma$ gives
functions of $z$ with simple poles at half-integers in some
left half-plane. Finally, a careful estimate of the size of
neighborhoods and remainders of the stationary phase method
shows that under the Diophantine condition the series over
$g$ converges as desired.
\end{proof}

\subsubsection{The odd case}

Now let $X$ be an odd-dimensional spin manifold with an isometric
action of $\Gamma$ such that all conditions from the preceding
subsection are satisfied. Let
$$
\sigma\in Mat(N,C^\infty(X)_\Gamma)
$$
be an invertible element. Consider the Toeplitz operator
$$
P\sigma P:P L^2(X,S\otimes \mathbb{C}^N)\longrightarrow P
L^2(X,S\otimes \mathbb{C}^N),
$$
where $P$ is the nonnegative spectral projection of the Dirac
operator $\mathcal{D}_X$.

\begin{theorem}
One has
\begin{equation}\label{local}
\ind P \sigma P=\sum_g\BL[\sum_{X_g}\;\;\int_{X_g}
{A(X_{g})}\Pf\left\{2i\sin\left(\frac\Omega{4\pi}+\frac{i\Theta}2\right)(NX_{g})\right\}^{-1}
\ch \sigma(g)\BR],
\end{equation}
where $\ch \sigma(g)$ is the coefficient of $g^*$ in the
Chern--Simons character of the invertible function $\sigma$.
\end{theorem}

\begin{proof}
This formula can be obtained from the corresponding even
formula obtained above by the standard trick
(multiplication of the manifold by $\mathbb{S}^1$;
cf.~\cite{MR813137}).
\end{proof}

\subsection{End of proof of the main theorem}

Now we can apply the results obtained in the preceding subsection
and prove the index formula.

Note that, assuming that $M$ is orientable, the cotangent
bundle of $M$ is a spin manifold; for the spin bundle one
can take the complexification of the lift of the exterior
form bundle from $M$ to $T^*M$. In this case, the Todd
operator is a Dirac operator.

Then the index formula \eqref{local} gives
\begin{multline}\label{local1}
\ind P \sigma(D)|_{S^*M}
P\\=\sum_g\BL[\sum_{(S^*M)_g}\;\;\int_{(S^*M)_g}\frac{A((S^*M)_g)
\ch_\Gamma \sigma(D)|_{S^*M}(g)}{\Pf
\left\{2i\sin\left(\frac\Omega{4\pi}+\frac{i\Theta}2\right)(N(S^*M)_{g})
\right\}} \BR].
\end{multline}

In our case, this gives the desired formula \eqref{fixedp}. Indeed,

1)  The fixed point set $(S^*M)_g$ is the cosphere bundle $S^*M_g$
of the corresponding fixed point set in $M$.

2) For the $A$-class of the fixed point set, one has (e.g., see
\cite{LaMi1})
$$
A(S^*M_g)=A(T^*M_g)=A(M_g)^2=\Td(TM_g\otimes \mathbb{C})
$$
(where equality holds not only for characteristic classes
but also for the corresponding differential forms).

3) Finally,
$$
\ch\lambda_{-1}(NM_g\otimes \mathbb{C})(g)= {\Pf
\left\{2i\sin\left(\frac\Omega{4\pi}
+\frac{i\Theta}2\right)(N(S^*M)_{g})\right\}}.
$$
This can be proved by more or less routine computations,
which will be given in the detailed version of the paper.

The proof of the index theorem is complete.


\end{document}